\documentclass[12pt,reqno]{amsart}  
\usepackage{enumerate}
\usepackage{url} 

\usepackage{ifthen} 

\usepackage[italian,english]{babel}
\usepackage{comment} 
\usepackage{fancyhdr}

\renewcommand{\qedsymbol}{$\square$}

\newcommand{\dx}{\mathrm{d}}
\newcommand{\eps}{\varepsilon}
\newcommand{\R}{\mathbb{R}}
\newcommand{\N}{\mathbb{N}}
\newcommand{\C}{\mathbb{C}}

\newcommand{\M}{\mathcal{M}}

\newcommand{\Odip}[2]{\mathcal{O}_{#1}\!\left(#2\right)\mathchoice{\!}{}{}{}}
\newcommand{\Odig}[1]{\mathcal{O}\Bigl(#1\Bigr)\mathchoice{\!}{}{}{}}

\newcommand{\Odim}[1]{\mathcal{O}\bigl(#1\bigr)}
\newcommand{\Odi}[1]{\Odip{}{#1}}

\newcommand{\odip}[2]{{o}_{#1}\!\left(#2\right)\mathchoice{\!}{}{}{}}
 
\newcommand{\odi}[1]{\odip{}{#1}}
 
\def \fine {{\hfill \qedsymbol}}
\def \CC {\mathbb C}

\newtheoremstyle{slantedtheorems}% name
{10pt}%      Space above
{6pt}%      Space below
{\slshape}%         Bo\dx y font
{}%         Indent amount (empty = no indent, \parindent = para indent)
{\bfseries}% Thm head font
{.}%        Punctuation after thm head
{.5em}%     Space after thm head: " " = normal interword space;
{\thmname{#1} \thmnumber{#2}\thmnote{ (#3)}}
 \theoremstyle{slantedtheorems}

  \newtheoremstyle{Boldtheorems}% name
  {10pt}%      Space above
  {6pt}%      Space below
  {\bfseries}%         Bo\dx y font
  {}%         Indent amount (empty = no indent, \parindent = para indent)
  {\bfseries}% Thm head font
  {.}%        Punctuation after thm head
  {.5em}%     Space after thm head: " " = normal interword space;
  {\thmname{#1} \thmnumber{#2}\thmnote{ (#3)}}
\theoremstyle{Boldtheorems}

\allowdisplaybreaks

\title[]{An extended pair-correlation conjecture \\ and primes in short intervals}
\author[]{A.~Languasco, A.~Perelli \lowercase{and} A.~Zaccagnini}
\numberwithin{equation}{section}

\begin{document}

\maketitle

\vskip.5cm
{\bf Abstract.} In this paper we extend the well-known investigations of Montgomery \cite{Mon/1973} and Goldston $\&$ Montgomery \cite{Go-Mo/1987}, concerning the pair-correlation function and its relations with the distribution of primes in short intervals, to a more general version of the pair-correlation function.



%-1-%%%%%%%%%%%%%%%%%%%%%%%%%%%%%%%%%%%%%%%%%%%%%%%
%%%%%%%%%%%%%%%%%%%%%%%%%%%%%%%%%%%%%%%%%%%%%%%%%
\section{Introduction}

\smallskip
This is a companion of our paper \cite{L-P-Z/tauPC}, where we study some problems on the distribution of primes assuming conjectural bounds for  the extended pair-correlation function
\[
F(X,T,\tau) = \sum_{-T\leq\gamma,\gamma'\leq T} X^{i(\gamma-\gamma')} w(\tau(\gamma-\gamma')),
\]
where $w(u) = 4/(4+u^2)$ and $0\leq \tau\leq 1$. We always assume the Riemann Hypothesis (RH) throughout the paper. We refer to the Introduction of \cite{L-P-Z/tauPC} for a discussion and motivations of the function $F(X,T,\tau)$. In this paper we extend to $F(X,T,\tau)$ the well-known results, under RH, of Montgomery \cite{Mon/1973} and of Goldston $\&$ Montgomery \cite{Go-Mo/1987}, respectively on the behavior of Montgomery's pair-correlation function
\[
F(X,T) = \sum_{-T\leq\gamma,\gamma'\leq T} X^{i(\gamma-\gamma')} w(\gamma-\gamma')
\]
and on the equivalence between Montgomery's conjecture for $F(X,T)$ and the asymptotic behavior of the mean-square of primes in short intervals. Motivation for the present work is to give some theoretical support to the assumptions we make in \cite{L-P-Z/tauPC}, and to refine the link between $F(X,T,\tau)$ and short averages of primes in short intervals.

\bigskip
We start with the extension to $F(X,T,\tau)$ of Montgomery's approach to $F(X,T)$ under RH. Throughout the paper, implicit constants in the $\ll, \mathcal{O}$ and $o$ symbols are either absolute or may depend on $\eps>0$ (if present in the formula), unless otherwise specified. Writing $f(x)=\infty(g(x))$ for $g(x)=o(f(x))$,
\begin{equation}
\label{1-1}
a(u,X,\tau)=  
\begin{cases}
(u/X)^{1/\tau} & \text{if}\ u \leq X \\ 
(X/u)^{1/\tau} & \text{if}\ u > X
\end{cases}
\quad \text{and} \quad S(X,\tau)=\sum_{n=1}^\infty \frac{\Lambda(n)^2}{n} a(n,X,\tau)^2
\end{equation}
we have

\medskip
%%%-TH1-%%%%%%%%%%%%%%%%%%%%
{\bf Theorem 1.} {\sl Let $X,T \geq 2$ and $\tau\in(0,1]$. Then $F(X,T,\tau)\geq 0$ and $F(X,T,\tau)=F(1/X,T,\tau)$.  Moreover, assuming RH, as $T\to \infty$ we have
\[
F(X,T,\tau)  = \frac{T}{\pi} \frac{S(X,\tau)}{\tau} \big(1+\odi{1}\big) + \frac{T \log^{2} T}{\pi \tau X^{2/\tau}}\big(1+\odi{1}\big) + \Odig{\frac{T  \log T(S(X,\tau))^{1/2}}{\tau X^{1/\tau}}} 
\]
uniformly for $\tau\geq 1/T$, provided that
\begin{equation} 
\label{1-2}
TS(X,\tau)=\infty(X)
\end{equation} 
and}
\begin{equation} 
\label{1-3}
\tau S(X,\tau)=\infty((\log^3 T)/T).
\end{equation} 

\medskip
{\bf Remark 1.} Since for $\tau=1$ we have $S(X,1)\sim \log X$ as $X\to\infty$ thanks to the prime number theorem, we see that Theorem 1 reduces to Montgomery's theorem \cite{Mon/1973} if $\tau=1$ and $X=T^\alpha$, $0< \alpha\leq1$. Indeed, in such a case the third error term is absorbed by the first two. However, in the general case we cannot compare such error terms since we don't have precise control on the size of $S(X,\tau)$; see below for more information on the function $S(X,\tau)$. \fine

\medskip
{\bf Remark 2.} We can relax \eqref{1-2} to the following group of milder conditions
\begin{equation}
\label{1-4}
T S(X, \tau) = \infty(X S(X, 2 \tau)), \qquad  T S(X, \tau) = \infty(\tau X), \qquad T S(X, \tau) = \infty(1/X^5).
\end{equation}
This apparently requires a slightly different treatment of the mean-square of $R_1(X,t,\tau)$, see \eqref{2-11} below, giving better results when $\tau\leq 1/\log X$. Note indeed that $X S(X, 2 \tau) \ll X$ in such a range thanks to \eqref{1-5} below, and that conditions \eqref{1-4} become $T=\infty(X)$ whenever the asymptotic formula \eqref{1-7} is available. We give a brief sketch of the argument at the end of the proof of Theorem 1. \fine

\bigskip
For $\tau=o(1)$, the function $S(X,\tau)$ depends on the distribution of primes in short intervals around $X$. Indeed, for $\eps>0$ and  $X\to \infty$, by the Brun-Titchmarsh inequality we have
\begin{equation}
\label{1-5}
S (X,\tau) \ll \tau \log X 
\end{equation}
uniformly for $X^{-1+\eps}\leq \tau \leq 1$. Moreover, given $\beta\in[0,1)$ and denoting by $K(\beta)$ the assertion
\begin{equation}
\label{1-6}
\psi(X+h)-\psi(X) \sim h \ \ \text{ uniformly for} \ \ X^{\beta+\eps} \leq h\leq X,
\end{equation}
if $K(\beta)$ holds true we have
\begin{equation}
\label{1-7}
S (X,\tau)  \sim \tau \log X  
\end{equation}
uniformly for $X^{-1+\beta+\eps}\leq \tau \leq 1$. The proof of \eqref{1-5} and \eqref{1-7} is standard, and we give a brief sketch at the end of the proof of Theorem 1.

\medskip
{\bf Remark 3.} It is well known that hypothesis $K(\beta)$ follows from RH for $\beta \geq 1/2$. Moreover, \eqref{1-5} and \eqref{1-7} show that for small $\tau$ there is a link between the behavior of $F(X,T,\tau)$ and the distribution of primes in short intervals. Such a link is already made explicit in \cite{L-P-Z/tauPC}, and is made more precise later in this paper. \fine

\medskip
From Theorem 1, Remark 2 and \eqref{1-7}  we obtain at once

\medskip
%%%-COR1-%%%%%%%%%%%%%
{\bf Corollary.} {\sl Let $X,T \geq 2$, and let  $0\leq \beta \leq 1/2$. Assuming that RH and $K(\beta)$ hold true, as $T\to \infty$ we have
\[
F(X,T,\tau) \sim \frac{T}{\pi} \log X
\]
uniformly for $T^\eps \leq X \leq T/\log T$ and} $\max  (X^{-1+\beta+\eps} ; T^{-1/2+\eps}) \leq \tau \leq 1$.

\bigskip
Lemma 2 of Heath-Brown $\&$ Goldston \cite{HB-Go/1984} gives an expression of $F(X,T,\tau)$ in terms of $F(X,T)$, which in our notation reads as
\begin{equation}
\label{1-8}
F(X,T,\tau) = \frac{F(X,T)}{\tau^2} + \frac{\tau^2 -1}{\tau^3} \int_0^{\infty} F(u,T) a(u,X,\tau)^2 \frac{\dx u}{u} \qquad \text{for} \ \tau\in(0,1].
\end{equation}
Hence we may plug in \eqref{1-8} a plausible quantitative version of Montgomery's pair-correlation conjecture in order to formulate a conjecture for $F(X,T,\tau)$ in the remaining range $X\geq T/\log T$. However, even plugging in \eqref{1-8} a rather sharp error term for $F(X,T)$, \eqref{1-8} allows to detect the uniform behavior of $F(X,T,\tau)$ only for $\tau$ quite close to 1. Hence \eqref{1-8} is useful to guess the main term of $F(X,T,\tau)$, but apparently it does not help much in the $\tau$-uniformity aspect.

\medskip
%%%-CONJECTURE-%%%%%%%%%
{\bf Conjecture.} {\sl Let $\eps>0$, $M>1$, $H=\min(X,T)$ and $K=\max(X,T)$. Then as $H\to\infty$}
\[
F(X,T,\tau) \sim \frac{T}{\pi} \log H  \ \  \text{uniformly for} \ \ K\leq H^M \ \ \text{and} \ \  H^{-1+\eps} \leq \tau \leq 1.
\]

\smallskip
{\bf Remark 4.} Note that, considering separately the two possibilities for $H$ and $K$, when $\tau=1$ the Conjecture coincides with Montgomery's conjecture plus, essentially, his result under RH. Moreover, part of the uniformity range of the  above asymptotic formula is already covered by the Corollary. Actually, when $H=X$ the Conjecture is supported by the Corollary, and condition $X^{-1+\eps}\leq \tau\leq 1$ is required in view of the erratic behavior of the sum $S(X,\tau)$ when $\tau$ is, roughly, of order $<1/X$. Moreover, lower bounds of type $\tau\geq (\log^AX)/X$ with any $A>1$, or even sharper, are suggested by Maier's \cite{Mai/1985} oscillation results, since \eqref{1-6} does not hold for all $X$ when $h\leq \log^AX$. When $H=T$, condition $T^{-1+\eps} \leq \tau \leq 1$ arises from a heuristic argument and, again, Maier's type results suggest that essentially wider $\tau$-ranges are forbidden; see Remark 5 in \cite{L-P-Z/tauPC}. We wish to thank Sandro Bettin and Adam Harper for exchanging ideas about the above mentioned heuristic argument and the $\tau$-ranges in the Conjecture. \fine

\medskip
Apart from the uniformity ranges, which are not immediately comparable, we see that the Conjecture is a refined version of Hypothesis $H(\eta)$ in \cite{L-P-Z/tauPC}, where we only assume an upper bound of type $F(X,T,\tau)\ll TX^\eps$. Note also that the sharper upper bound $F(X,T,\tau)\ll T\log X$ follows from Theorem 1 and \eqref{1-5} in certain ranges, which however are disjoint from those needed in \cite{L-P-Z/tauPC}. 

\bigskip
In \cite{L-P-Z/tauPC} we deduce from $H(\eta)$ sharp bounds for short mean-square averages of primes in short intervals; here we show that the Conjecture is actually equivalent to the asymptotic behavior of such mean-squares, at least in a short range of $\tau$ close to 1. More precisely, the Conjecture determines the behavior of
\[
J(X,\tau,\theta)=\int_X^{X(1+\tau)}\bigl(\psi(x+\theta x) - \psi(x) -\theta x  \bigr)^2  \dx x \hskip1.5cm 0\leq \theta,\tau\leq1
\]
in a relatively large range of $\tau$, but the opposite implication requires much stronger limitations. This, as well as the basic condition $T^{-1/2+\eps} \leq \tau\leq 1$ in the Corollary, is due to the use of trivial bounds for $F(X,T,\tau)$ and for the related function $\Phi(X, t, \tau)$ in \eqref{2-1}, which for small $\tau$ are much worse than the expected order. Prototypical examples are Lemma 3 and  the error term in  \eqref{2-10}. 

\medskip
We refer to \cite{L-P-Z/tauPC} for the classical results on $J(X,\theta) = J(X,1,\theta)$ under RH. The implication from the Conjecture to $J(X,\tau,\theta)$ is given by

\medskip
%-TH2-%%%%%%%%%%%%
{\bf Theorem 2.} {\sl Assume RH and the Conjecture, and let $\eps>0$. Then as $X\to\infty$
\begin{equation}
\label{1-9} 
J(X, \tau, \theta) \sim \Bigl(1+\frac{\tau}{2}\Bigr)\, \tau  \theta X^{2} \log(1/\theta)
\end{equation}
uniformly for $1/X \leq \theta \leq X^{-\eps}$ and} $\theta^{1/2-\eps}\leq \tau \leq 1$.

\medskip
Actually, in the proof of Theorem 2 we use the Conjecture only in the ranges $T/\log T \leq X \leq T^M$ and $T^{-1/2+\eps} \leq \tau\leq 1$. In the opposite direction we have the following weaker result. We omit its proof, which follows the lines of the corresponding result in  Goldston $\&$ Montgomery \cite{Go-Mo/1987}, suitably modified as for Theorem 2. Assume RH and let $\eps>0$ and $X\to\infty$. If \eqref{1-9} holds uniformly for  $1/X \leq \theta \leq \tau^{-1} X^{-\eps}$ and $\max (X^{-\eps}; \theta^{1/2-\eps})\leq \tau \leq 1$, then
\[
F(X, T, \tau)  \sim \frac{T}{\pi}\log T 
\]
uniformly for $X^{\eps} \leq T \leq X$ and $T^{-1/2+\eps}\leq \tau \leq 1$. Note that, in both cases, for $\tau=1$ we get back the results in \cite{Go-Mo/1987}.

\medskip
{\bf Remark 5.} By the arguments of Saffari $\&$ Vaughan \cite{Sa-Va/1977} and Goldston $\&$ Montgomery \cite{Go-Mo/1987} one can show that, as $X\to\infty$, the asymptotic formula \eqref{1-9}, uniformly for $X^{-1+\eps} \leq \theta \leq X^{-\eps}$ and $\theta\leq \tau \leq 1$, is equivalent to
\[
\int_X^{X+Y}\bigl(\psi(x+h) - \psi(x) - h  \bigr)^2 \dx x \sim Yh \log(X/h),
\]
uniformly for $X^\eps \leq h \leq X^{1-\eps}$ and $h \leq Y \leq X$. \fine

\medskip
We conclude observing that it would be helpful having numerical evidence and/or Random Matrix Theory heuristics supporting the above Conjecture. However, apparently both are not so easy to obtain, at least with the present form of the Conjecture.

\medskip
{\bf Acknowledgements.} 
%%%%%%%%%%%%%
We thank Sandro Bettin, who read a previous version of the paper, pointing out inaccuracies and suggesting improvements at several places. We also thank the referee, who pointed out further inaccuracies and suggested improvements in the presentation. This research was partially supported by the grant PRIN2010-11 {\sl Arithmetic Algebraic Geometry and Number Theory}.

\bigskip
%-2-%%%%%%%%%%%%%%%%%%%%%%%%%%%%%%%%%%%%%%%%%%%%%%%
%%%%%%%%%%%%%%%%%%%%%%%%%%%%%%%%%%%%%%%%%%%%%%%%%
\section{Proof of Theorem 1}

\smallskip
We recall that the implicit constants in the $\ll, \mathcal{O}$ and $o$ symbols are either absolute or may depend on $\eps>0$, unless otherwise specified. The first assertion of Theorem 1 follows  from
\[ 
\int_{-\infty}^{+\infty} \vert \sum\limits_{-T \le \gamma  \le  T} X^{i\gamma} e^{i\gamma v}\vert^2 e^{-2 \vert v \vert /\tau} \dx v = \tau F(X,T,\tau),
\]
see e.g. Lemma 2 of \cite{L-P-Z/tauPC}, while the second is trivial. Since for the main assertion of Theorem 1 we follow the proof in Montgomery \cite{Mon/1973}, we shall be sketchy, reporting the structure of the proof and giving details only when we have to keep track of the uniformity in $\sigma$, which will later be related with $\tau$. The starting point of the proof is a slightly modified form of the Lemma on page 185 of \cite{Mon/1973}, see Lemma 4 below. We need three auxiliary lemmas.

\medskip
%%%%%%%%%%%%%%%%%%%%%
{\bf Lemma 1.} {\sl Let $D=\{ s \in \C: \Re(s) < 1, |\Im(s)| < 1\}$ and $s \in \C \setminus D$. Then}
\[
\frac{\Gamma'}{\Gamma}(s) = \log s + \Odi{1}.
\]

\medskip
{\it Proof.} The starting point is the following integral representation (see formula $\Gamma15$ on page 428 of Lang \cite{Lan/1999}), which holds for every $s \in \C \setminus D$:
\[
  \frac{\Gamma'}{\Gamma}(s)  = \log s - \frac1{2 s} +  2 \int_0^{+\infty} \frac{P_2(y)}{(s + y)^3} \dx y,
\]
where $P_2(y)$ is the periodic function of period $1$ which equals $\frac12 (y^2 - y)$ for $y \in [0, 1]$. Lemma 1 follows then by simple estimates of the integral, considering separately the two cases $|t|\geq 1$ and $|t|\leq 1$, $\sigma\geq 1$. \fine

\medskip
%%%%%%%%%%%%%%%%%%%%
{\bf Lemma 2.} {\sl For $a,b$ with $1 \le a \le b$, $\tau \in (0, 1]$ and $\kappa \ge 0$ we have
\[
\sum_{n \in [a,b]} \frac{\log^\kappa n}{1 + \tau^2 n^2} \ll
\begin{cases}
(b-a + 1) \log^\kappa b \ \ &\text{if $a \le b \le \tau^{-1}$} \\ 
\tau^{-1} \log^\kappa(2 / \tau) &\text{if $a \le \tau^{-1} \le b$} \\
(\tau^2 a)^{-1} \log^\kappa a &\text{if $\tau^{-1} \le a \le b $,}
  \end{cases}
\]
where the implicit constant depends at most on $\kappa$.}

\medskip
{\it Proof}. The first inequality is trivial. To prove the second inequality we split the interval $[a,b]$ into $[a, \tau^{-1}] \cup [\tau^{-1}, b]$. As before, the sum over $[a, \tau^{-1}]$ is $\ll \tau^{-1} \log^{\kappa}(2 / \tau)$, while
\[
  \sum_{n \in [\tau^{-1}, b]}  \frac{\log^\kappa n}{1 + \tau^2 n^2}  \ll  \tau^{-2} \int_{\tau^{-1}}^b \frac{\log^\kappa t}{t^2} \dx t  \ll  \tau^{-1} \log^\kappa(2/\tau).
\]
The third case is similar, with the interval $[a,b]$ in place of $[\tau^{-1}, b]$. \fine

\medskip
We shall use Lemma 2 only with $\kappa=0$ or $\kappa=1$. Let
\begin{equation}
\label{2-1} 
\Phi(X, t, \tau)=\Bigl \vert  \sum_{\gamma}\frac{X^{i\gamma}}{1+\tau^{2}(t-\gamma)^{2}}\Bigr \vert.
\end{equation}

\smallskip
%%%%%%%%%%%%%
{\bf Lemma 3.} {\sl Let $X\geq 2$, $t\in \R$ and $\tau \in (0,1]$, and write $M=\max  (\vert t \vert +2; 2/ \tau)$. Then}
\[
\Phi(X, t, \tau) \leq \Phi(1, t, \tau) \ll \frac{1}{\tau}\log M.
\]

\medskip
{\it Proof.} By the Riemann-von Mangoldt formula and Lemma 2 we have
\[
\Phi(1, t, \tau)  \ll \sum_{n\geq 1} \frac{\log (n+ \vert t \vert)}{1+\tau^{2} n^2} \ll \log M \sum_{n\leq M}
\frac{1}{1+\tau^{2} n^2} + \sum_{n > M} \frac{\log n}{1+\tau^{2} n^2} \ll \frac{\log M}{\tau},
\]
and  Lemma 3 follows. \fine

\medskip
%%%%%%%%%%%%%%%%%%%%%%%%%
{\bf Lemma 4} {\sl Assume RH and let $X \geq 2$. Then, uniformly for $\sigma>1$ and $|t|\geq 1$, we have}
\begin{align*}
(2\sigma-1)  \sum_{\gamma} &\frac{X^{i\gamma}}{(\sigma -1/2)^{2}+(t-\gamma)^{2}} = -X^{-1/2}\Bigl(
\sum_{n\leq X} \Lambda(n) \Bigl( \frac{X}{n} \Bigr)^{1-\sigma+it}  + \sum_{n > X} \Lambda(n) \Bigl( \frac{X}{n} \Bigr)^{\sigma+it} \Bigr)\\
&+ X^{1/2-\sigma+it} \Bigl( \frac{\zeta'}{\zeta}(\sigma-it)+\frac{1}{4} \log (\sigma^{2} + t^{2} )+\frac{1}{4} \log ((1-\sigma)^{2} + t^{2})\Bigr) \\
&+ \frac{X^{1/2}}{\sigma-1+it} + \frac{X^{1/2}}{\sigma-it} +\Odi{\frac{1}{X^{5/2}\vert t \vert}} + \Odi{X^{1/2-\sigma} }.
\end{align*}

\medskip
{\it Proof.} We follow the proof of the Lemma in \cite{Mon/1973} till equation (22), which in our notation reads as
\begin{align}
\notag
(2\sigma-1) \sum_{\gamma}
&\frac{X^{i\gamma}}{(\sigma -1/2)^{2}+(t-\gamma)^{2}}= -X^{-1/2} \Bigl(\sum_{n\leq X} \Lambda(n) \Bigl( \frac{X}{n} \Bigr)^{1-\sigma+it} +\sum_{n > X} \Lambda(n) \Bigl( \frac{X}{n} \Bigr)^{\sigma+it} \Bigr)\\
\notag
& -\frac{\zeta'}{\zeta}(1-\sigma+it) X^{1/2-\sigma+it} +\frac{X^{1/2}(2\sigma-1)}{(\sigma-1+it)(\sigma-it)}\\
\label{2-2}
&-X^{-1/2}(2\sigma-1)\sum_{n=1}^{+\infty}\frac{X^{-2n}}{(\sigma-1-it-2n)(\sigma+it+2n)}.
\end{align}
In order to treat the second term on the r.h.s. of \eqref{2-2} we consider the logarithmic derivative of the functional equation of $\zeta(s)$ and use Lemma 1 to compute the resulting $\Gamma'/\Gamma$-terms. Since the argument of $s$ is bounded we have
\[
-\frac{\zeta'}{\zeta}(1-\sigma+it) = \frac{\zeta'}{\zeta}(\sigma-it) +\frac{1}{4} \log (\sigma^{2} + t^{2} ) +\frac{1}{4} \log ((1-\sigma)^{2} + t^{2} ) + \Odi{1},
\]
which we insert in \eqref{2-2}. By trivial estimates we have
\[
\sum_{n \ge 1} \frac{X^{-2 n}}{(\sigma - 1 - i t - 2 n)(\sigma + i t + 2 n)}  \ll  \sum_{n \ge 1}\frac{X^{-2n}}{|t|((\sigma+2n)^2+t^2)^{1/2}} \ll \frac{1}{|t|\sigma X^2},
\]
hence the last term in \eqref{2-2} is $\ll \frac{2 \sigma - 1}{\sigma} X^{- 5 /2} \vert t \vert^{-1}$, and the proof of Lemma 4 is complete. \fine

\bigskip
{\bf Proof of Theorem 1}. In the estimates below we tacitly assume condition $\tau\geq 1/T$. We apply Lemma 4 with $\sigma=\sigma_0 = 1/2+1/\tau$ and $|t| \geq 1$, and write the resulting identity as 
\begin{equation}
\label{2-3}
L(X,t,\tau) = R(X,t,\tau).
\end{equation}
We start by computing the $L^2$-norm of  $L(X,t,\tau)$ over $\mathcal{I}(T) =[-T,-1]\cup [1,T]$. Recalling \eqref{2-1} we have
\[
\int_{\mathcal{I}(T)} \vert L(X,t,\tau) \vert^{2} \dx t = 4 \tau^{2} \int_{\mathcal{I}(T)} \Phi(X, t, \tau)^2 \dx t
\]
and by Lemma 3 we get
\[
\int_{-1}^{1}\Phi(X, t, \tau)^2 \dx t \ll\frac{\log^2 T} {\tau^{2}},
\]
therefore
\begin{equation}
\label{2-4}
\int_{\mathcal{I}(T)}\vert L(X,t,\tau)\vert^{2} \dx t = \tau {\mathcal J}(X,T,\tau) +\Odi{\log^2 T}
\end{equation}
with
\begin{equation}
\label{2-5}
{\mathcal J}(X,T,\tau) = 4 \tau\int_{-T}^{T} \Phi(X, t, \tau)^2 \dx t.
\end{equation}
Next we link ${\mathcal J}(X,T,\tau)$ to $F(X,T,\tau)$. Let $|t|\leq T$. By the Riemann-von Mangoldt formula and Lemma 2 we have
\begin{align*}
\sum_{\vert \gamma \vert >T}\frac{1}{1+\tau^{2}(t-\gamma)^{2}}
&\ll \frac{\log T}{\tau^{2}} \frac{1}{T-t+1}
\end{align*}
and thanks to Lemma 3 we get 
\begin{equation}
\label{2-6}
\int_{-T}^{T} \sum_{\gamma, \gamma' \colon \vert \gamma \vert>T} \frac{\dx t}{(1+\tau^{2}(t-\gamma)^{2})(1+\tau^{2}(t-\gamma')^{2})} \ll \frac{\log^{3}T}{\tau^{3}}.
\end{equation}
For $\vert t \vert > T$, a similar computation based on Lemma 2 shows that 
\begin{equation}
\label{2-7}
\int_{\{\vert t \vert > T\}} \sum_{-T\leq \gamma, \gamma' \leq T}  \frac{\dx t}{[1+\tau^{2}(t-\gamma)^{2}][1+\tau^{2}(t-\gamma')^{2}]} \ll \frac{\log^{2}T}{\tau^{3}},
\end{equation}
hence inserting \eqref{2-6} and \eqref{2-7} in \eqref{2-4} we have
\begin{equation}
\label{2-8}
\int_{\mathcal{I}(T)} \vert L(X,t,\tau) \vert^{2} \dx t = 4 \tau^{2} \int_{-\infty}^{+\infty} \Bigl \vert \sum_{-T\leq \gamma \leq T}
\frac{X^{i\gamma}}{1+\tau^{2}(t-\gamma)^{2}} \Bigr \vert^{2} \dx t + \Odig{\frac{\log^{3}T}{\tau}}.
\end{equation}
Computing residues as on p.188 of \cite{Mon/1973} we see that
\[
\int_{-\infty}^{+\infty} \frac{\dx t}{(1+\tau^{2}(t-\gamma)^{2})(1+\tau^{2}(t-\gamma')^{2})} 
= \frac{\pi}{2\tau} \frac{4}{4+\tau^{2}(\gamma - \gamma')^{2}} =  \frac{\pi}{2\tau} w(\tau(\gamma-\gamma')),
\]
hence from \eqref{2-8} and the definition of $F(X,T,\tau)$ we finally obtain
\begin{equation} 
\label{2-9}
\int_{\mathcal{I}(T)} \vert L(X,t,\tau) \vert^{2} \dx t = 2 \pi \tau F(X,T,\tau) + \Odig{\frac{\log^{3}T}{\tau}}.
\end{equation}
For future reference we remark that, in view of  \eqref{2-4}, \eqref{2-9} may be expressed as
\begin{equation}
\label{2-10}
2 \pi  F(X,T,\tau) = {\mathcal J}(X,T,\tau) + \Odig{\frac{\log^{3}T}{\tau^2}}.
\end{equation}

\medskip
We now turn to $R(X,t,\tau)$. We first write
\begin{equation}
\label{2-11}
R(X,t,\tau) = \sum_{j=1}^{5}R_{j}(X,t,\tau)
\end{equation}
where, in view of Lemma 4 and recalling that $\sigma_0=1/2+1/\tau\geq 3/2$ and $\vert t \vert \geq 1$,
\begin{align*}
R_{1}(X,t,\tau) & =  -X^{-1/2} \Bigl(\sum_{n\leq X} \Lambda(n) \Bigl( \frac{X}{n} \Bigr)^{1-\sigma_0+it} + 
\sum_{n > X} \Lambda(n) \Bigl( \frac{X}{n} \Bigr)^{\sigma_0+it}\Bigr), \\
R_{2}(X,t,\tau) & = X^{1/2-\sigma_0+it}   \frac{\zeta'}{\zeta}(\sigma_0-it) + \Odi{X^{1/2-\sigma_0}} \ll X^{1/2-\sigma_0}, \\
R_{3}(X,t,\tau) & = \frac{1}{4} X^{1/2-\sigma_0+it} \Bigl(\log (\sigma_0^{2} + t^{2} ) + \log ((1-\sigma_0)^{2} + t^{2} )\Bigr), \\
R_{4}(X,t,\tau) & =\frac{X^{1/2}}{\sigma_0-1+it} + \frac{X^{1/2}}{\sigma_0-it}, \qquad R_{5}(X,t,\tau) \ll \frac{1}{ X^{5/2}\vert t \vert}.
\end{align*}
Next we compute
\begin{equation}
\label{2-12}
M_j(X,T,\tau)= \int_{\mathcal{I}(T)} \vert R_{j}(X,t,\tau) \vert^{2} \dx t,
\end{equation}
but first we recall a mean-value theorem for exponential sums (Lemma 6 of Goldston $\&$ Montgomery \cite{Go-Mo/1987}) and a sieve upper bound for $k$-twin primes (see Theorem 3.11 of Halberstam $\&$ Richert \cite{Ha-Ri/1974} for the sieve bound and Lemma 17.4 of Montgomery \cite{Mon/1971} for the summation of the singular series). Write $e(x)=e^{2\pi ix}$ and let $\M$ be a countable set of real numbers, $k\in\N$ and
\[
S(t)=\sum_{\mu \in \M} c(\mu) e(\mu t), \hskip1.5cm Z(U; k) = \sum_{n\leq U} \sum_{\substack{m\leq U \\  |n-m| = k}} \Lambda(n) \Lambda(m).
\]

\smallskip
{\bf Lemma 5.} {\sl Let $T\geq 1$, $1/(2T) \leq \delta \leq 1/2$ and $S(t)$ be as above with $c(\mu)\in \R$ and $\sum_{\mu \in \M} \vert c(\mu) \vert $ convergent. Then}
\[
\int_{-T}^T \vert S(t) \vert ^2 \dx t = \Bigl(2T+\Odig{\frac{1}{\delta}}\Bigr) \sum_{\mu \in \M} \vert c(\mu) \vert ^2 + \Odig{T
\sum_{\substack{\mu,\nu \in \M \\ 0< \vert \mu - \nu \vert <\delta}}  \vert c(\mu)c(\nu) \vert }.
\]

\smallskip
{\bf Lemma 6.} {\sl Let  $U,V\geq 1 $. Then} $\sum_{k\leq V}Z(U;k) \ll UV$.

\medskip
By \eqref{1-1} and Lemma 5 we obtain
\begin{equation}
\label{2-13}
\int_{-T}^{T} \vert  R_{1}(X,t,\tau) \vert^{2} \dx t =  2 T S(X,\tau) + \Odig{ \delta^{-1}S(X,\tau) + E_{1} + E_{2}},
\end{equation}
where
\begin{align*}
E_{1} &= TX^{-1/\tau}\sum_{n\leq X}\sum_{0< \vert \log(n/m) \vert <2\pi\delta}\Lambda(n)\Lambda(m) n^{1/\tau-1/2} \frac{a(m,X,\tau)}{m^{1/2}}, \\
E_{2} &= TX^{1/\tau}\sum_{r =0}^{+\infty}\sum_{n = 2^r X}^{2^{r+1}X}\sum_{0< \vert \log(n/m) \vert <2\pi\delta}
\Lambda(n) \Lambda(m) n^{-1/\tau-1/2} \frac{a(m,X,\tau)}{m^{1/2}}.
\end{align*}
Hence by Lemma 6 and a standard dissection argument we get the bounds
\begin{equation}
\label{2-14}
\begin{split}
E_{1} &\ll TX^{-1}\sum_{k \leq 10^4\delta X}Z(10^4X;k) \ll T\delta X \\
E_{2} &\ll TX^{-1} \sum_{r =0}^{+\infty} 2^{(-2/\tau-1)r} \sum_{k \leq 10^4  2^{r+1}\delta X} Z(10^4 2^{r+1} X;k) \ll
T\delta X.
\end{split}
\end{equation}
Choose $2\delta=(TX)^{-1/2}(S(X,\tau))^{1/2}$. Note that $S(X,\tau)$ is increasing in $\tau$, hence from standard prime number theory we have $S(X,\tau)\leq S(X,1) \ll \log X$ and therefore $\delta \leq 1/2$. Moreover, from \eqref{1-2} we deduce that $\delta\geq 1/(2T)$ since $T\to\infty$. In the rest of the proof we tacitly use \eqref{1-2} and \eqref{1-3}. From $1/(2T) \leq \delta \leq 1/2$, inserting \eqref{2-14} into \eqref{2-13} we get
\begin{equation}
\label{2-15}
\int_{-T}^{T} \vert R_{1}(X,t,\tau)  \vert^{2} \dx t = 2 T S(X,\tau) + \Odi{(T X S(X,\tau))^{1/2} }.
\end{equation}
Since a direct estimate implies that  
\[
  \int_{-1}^1 \vert R_{1}(X,t,\tau) \vert^2 \dx t  \ll X,
\]
in view of \eqref{2-15} we have
\begin{equation}
\label{2-16}
M_{1}(X,T,\tau)=2 T S(X,\tau)+\Odi{(T X S(X,\tau))^{1/2} }.
\end{equation}
Coming to the terms with $j\geq 2$, simple computations show that
\begin{equation}
\label{2-17}
M_{2}(X,T,\tau)\ll X^{-2/\tau} T, \hskip1cm M_{3}(X,T,\tau) =  2 X^{-2/\tau} T \log^{2} T + \Odi{X^{-2/\tau} T \log T},
\end{equation}
while, thanks to \eqref{1-2},
\begin{equation}
\label{2-18}
M_{4}(X,T,\tau) \ll \tau X \int_{-\tau T}^{\tau T} \frac{\dx t}{1+t^{2}}  \ll \tau X = \odi{T S(X, \tau)}
\end{equation}
and
\begin{equation}
\label{2-19} 
M_{5}(X,T,\tau) \ll \frac{1}{X^{5}} \int_{1}^{T} \frac{\dx t}{t^{2}} \ll \frac{1}{X^{5}} = \odi{TS(X,\tau)}.
\end{equation}

\medskip
Finally, by \eqref{2-11}, \eqref{2-12} and the Cauchy-Schwarz inequality we have
\[
\int_{\mathcal{I}(T)} \vert R(X,t,\tau) \vert ^{2} \dx t =\sum_{j=1}^{5} M_{j}(X,T,\tau)+\Odig{\sum_{j=1}^{5}
\sum_{\substack{k=1 \\ j\neq k}}^{5} \big(M_j(X,T,\tau)M_k(X,T,\tau)\big)^{1/2}},
\]
hence from \eqref{2-16}-\eqref{2-19} we get
\begin{align} 
\notag
\int_{\mathcal{I}(T)} \vert R(X,t,\tau) \vert ^{2} \dx t &= 2 T  S(X,\tau) + \frac{2 T  \log^{2} T}{X^{2/\tau}} + \Odig{\frac{T  \log T( S(X,\tau) )^{1/2}}{X^{1/\tau}}} \\
\label{2-20}
&+\Odig{\frac{T  \log T}{X^{2/\tau}} +(T  X  S(X,\tau) )^{1/2} } +  \odi{T S(X, \tau)}.
\end{align}
In view of \eqref{2-3}, \eqref{2-9} and \eqref{2-20}, dividing by $2 \pi \tau$ we obtain
\begin{align}
\notag
F(X,T,\tau) &= \frac{T}{\pi\tau} S(X,\tau) + \frac{T \log^{2} T}{\pi \tau X^{2/\tau}} + \Odig{\frac{T  \log T(S(X,\tau))^{1/2}}{\tau X^{1/\tau}}} \\
\notag
&+ \Odig{\frac{T  \log T}{\tau X^{2/\tau}} + \frac{(T  X S(X,\tau))^{1/2}}{\tau}  + \frac{\log^{3}T}{\tau^{2}}} +  \odi{\frac{T}{\tau} S(X, \tau)}.
\end{align}
Note that the first error term is the geometric mean of the two main terms. Hence it is negligible when one of the main terms dominates, but not when they have the same size. This ends the proof of Theorem 1, since the other error terms are of the required size provided \eqref{1-2} and \eqref{1-3} hold. \fine

\bigskip
For the alternative treatment mentioned in Remark 2 we write
\[
\widetilde{S}(X, \tau) = \sum_{n \ge 1} \Lambda(n)^2 \, a(n, X, \tau)^2,
\]
thus Corollary 3 of Montgomery \& Vaughan \cite{Mon-Vau/1974} yields
\[
\int_{-T}^T \vert R_1(X, t, \tau) \vert^2 \, \dx t = \sum_{n \ge 1} \bigl( 2 T + \Odi{n} \bigr)  \frac{\Lambda(n)^2}n \, a(n, X,\tau)^2 = 2 T S(X, \tau)+  \Odim{\widetilde{S}(X, \tau)}.
\]
Moreover, a computation shows that $\widetilde{S}(X, \tau) \le X S(X,2 \tau)$. Hence we may apply the above formula (also for $T=1$), thus getting
\[
M_{1}(X,T,\tau)=2 T S(X,\tau)+\Odi{X S(X,2\tau)}.
\]
The treatment of $M_j(X,T,\tau)$ for $j=2,\dots,5$ is as before, and conditions \eqref{1-4} in Remark 2 lead to the bounds \eqref{2-17}--\eqref{2-19}. The assertion in Remark 2 now follows arguing as before. \fine

\bigskip
We now turn to a brief sketch of \eqref{1-5} and \eqref{1-7}. Given a parameter $1\leq H\leq X$ we split the range of summation in $S(X,\tau)$ as $[1,X-H] \cup (X-H,X+H] \cup (X+H,\infty)$ and denote by $S_1,S_2,S_3$ the corresponding subsums. Choosing $H = \tau X \in [X^\eps, X]$, by the trivial estimate for $\Lambda(n)$ and the Brun-Titchmarsh inequality  we obtain
\[
S_2 \ll \frac{\log X}{X} (\psi(X+H)-\psi(X-H)) \ll \frac{H \log^{2} X}{X \log H} \ll \tau \log X.
\]
The same bound can be obtained for the subsums $S_1$ and $S_3$, by a further splitting-up argument into intervals of length at most $H$ and then proceeding similarly, thus getting \eqref{1-5}. Assume now hypothesis $K(\beta)$ and choose $H=h$. Trivial estimates then give
\begin{equation}
\label{2-21}
S_{1} + S_3 = \Odig{\tau\log^2X\Bigl(1-\frac{H}{X}\Bigr)^{2/\tau}},
\end{equation}
while by partial summation we get
\begin{equation}
\label{2-22}
S_{2} = \tau \log X \Bigl(1-\Bigl(1-\frac{H}{X}\Bigr)^{2/\tau}\Bigr) + \Odig{\frac{H^2}{X^2} + \frac{\log^3X}{X}(1+HX^{-1/2})} + \odi{\tau \log X }.
\end{equation}
Hence, letting, say, $H=\tau X (\log X)^{1/3}$ for $X^{-1+\beta+\eps}\leq \tau \leq (\log X)^{-1/2}$ and $H= X (1- 1/(\log X)^{\tau})$ otherwise, we have that $h=H$ is consistent with \eqref{1-6}, and \eqref{1-7} follows from \eqref{2-21} and \eqref{2-22}. \fine

\medskip
Similar computations show that the behavior of $S(X,\tau)$ becomes erratic (depending essentially on the prime-power closest to $X$) when $\tau$ is, roughly, of order $<1/X$.

\bigskip
%-SECT3-%%%%%%%%%%%%%%%%%%%%%%%%%%%%%%%%%%%
%%%%%%%%%%%%%%%%%%%%%%%%%%%%%%%%%%%%%%%%
\section{Proof of Theorem 2}

\smallskip
We start with several lemmas, which are a $\tau$-uniform version of the Goldston $\&$ Montgomery \cite{Go-Mo/1987} abelian-tauberian method. We skip several details, referring instead to \cite{Go-Mo/1987} or \cite{L-P-Z/2012}.

\medskip
%-L7-%%%%%%%%%%%%%%%
{\bf Lemma 7.} {\sl Let $f(y) \ge 0$ be a continuous function of $y\in \R$ and $\tau\in (0,1]$. Suppose that
\[
 I(Y)  =  \int_{-\infty}^{\infty} e^{-2 |y|} f(Y + y) \dx y = 1 + \eps(Y)
\]
with $|\eps(Y)| \leq 1/2$. Then, for every $a,b$ with $0\leq a < b\leq 1$,
\[
\int_a^b e^{2 y \tau} f(Y + y) \dx y=  \Bigl(\int_a^b e^{2 y \tau} \dx y \Bigr) \Bigl(1+ \Odi{ \overline{\eps}(Y)^{1/2}}\Bigr)
\] 
as $Y\to\infty$, where} $\overline{\eps}(Y) =   \sup_{t \in [a - 1, b + 1]}|\eps(Y + t)|$.

\medskip
{\it Proof.} We follow the proof of Lemma 1 of \cite{L-P-Z/2012} until eq.~(3.8). So  we may write
\[
 \int_a^y f(Y + t) \dx t = y - a + \Odi{\overline\eps(Y)^{1/2}}
\]
and hence a partial integration argument shows that
\[
\int_a^b (f(Y + y) -1) e^{2y \tau} \dx y  \ll \overline\eps(Y)^{1/2} \int_a^b e^{2y \tau} \dx y  
\]
since $a,b$ and $\tau$ are bounded, and Lemma 7 follows.    \fine 

\medskip
Recalling \eqref{2-1}, in the following lemma we link $F(X,T,\tau)$ with
\[
I(X,\tau,\kappa) = 4\tau \int_{0}^{\infty} \Bigl(\frac{\sin \kappa t}{t}\Bigr)^{2} \Phi(X, t, \tau)^2 \dx t
\]
by means of relation \eqref{2-10}, connecting $F(X,T,\tau)$ with the integral $\mathcal{J}(X,T,\tau)$ defined in \eqref{2-5}. Such a lemma is a $\tau$-uniform version of Lemma 2 of  \cite{Go-Mo/1987}, the main difference being that we avoid the use of individual bounds for  $\Phi(X, t, \tau)$, which are weak for small $\tau$. We have

%-L8-%%%%%%%%%%%%%%%%%%%%%
\medskip
{\bf Lemma 8.} {\sl Assume the Conjecture and let $\eps>0$. Then as $X\to\infty$
\[ 
I(X,\tau,\kappa) \sim \frac{\pi}{2} \kappa \log (1/\kappa)
\]
uniformly for $1/X \leq \kappa \leq X^{-\eps}$ and $\kappa^{1/2-\eps}\leq \tau\leq 1$.} 

\medskip
{\it Proof.} For simplicity we write $\Phi(t) = 4\tau \Phi(X, t, \tau)^2$. Since $\Phi(t)$ is an even function we may restrict to $t\geq 0$. We write $r(t)= (1/2){\mathcal J}(X,t,\tau)-t\log t$ and observe that $r'(t) = \Phi(t) -\log t-1$ and also $\Phi(t)=\log (1/\kappa)+1+\log(\kappa t) + r'(t)$. Hence, thanks to formulae 3.821.9 and 4.423.3 of Gradshteyn \& Ryzhik \cite{GradshteynR2007}, we have
\begin{equation}
\label{3-1}
I(X,\tau,\kappa)  =\frac{\pi}{2} \kappa  \log(1/\kappa) + \Odi{\kappa} + \int_{0}^{\infty} \Bigl(\frac{\sin \kappa t}{t}\Bigr)^{2} r'(t) \dx t.
\end{equation}
In order to estimate the integral in \eqref{3-1} we split $[0,\infty)$ as $[0,U] \cup [U,V] \cup [V,\infty)$, where $0<U<V$ are chosen below, and denote by $I_1,I_2,I_3$ the resulting integrals. Note that we may apply our Conjecture to $F(X,t,\tau)$ for every $U\leq t\leq V$ provided $\tau\geq U^{-1/2+\eps}$ and $X^\eps\leq U,V\leq X^A$ for some $A>0$. Since $\eps>0$ is arbitrary, the choice
\begin{equation}
\label{3-2}
U = \big(\kappa\log(1/\kappa)\big)^{-1},  \hskip1.5cm V = X^2
%(\tau\kappa)^{-1}\log^{2} (1/\kappa)
\end{equation}
satisfies the above inequalities provided $\kappa$ and $\tau$ belong to the ranges in the statement of the lemma. Therefore, since $U<X$ and the error term in \eqref{2-10} is $o(t\log t)$ for $t$ and $\tau$ as above, we have that if $f(X)\to\infty$ slowly enough then
\begin{equation}
\label{3-3}
r(t) = 
\begin{cases}
 o(t\log t) \  &\text{if} \ U\leq t \leq Xf(X) \\
 \Odi{t\log t}  \   &\text{if} \ Xf(X) \leq t \leq V.
\end{cases}
\end{equation}

\smallskip
In what follows we shall repeatedly use, without further mention, the bounds 
\[
\sin (\kappa t)/t \ll \min(\kappa, 1/ \vert t \vert) \hskip1.5cm \frac{\partial}{\partial t}(\sin (\kappa t)/t)^{2} \ll \kappa \vert t\vert^{-1} \min (\kappa;\vert t\vert^{-1})
\]
as well as the choice \eqref{3-2} and relation \eqref{2-10}. Using $\sin (\kappa t)/t \ll \kappa$ we have
\begin{equation}
\label{3-4}
\begin{split}
I_1 &\ll \kappa^2 \int_{0}^{U}  \big(\Phi(t) + |\log t| + 1\big) \dx t \\
&\ll \kappa^2\Big(F(X,U,\tau) + \frac{\log^3U}{\tau^2}  + \int_0^1|\log t|\dx t + U\log U\Big) \ll \kappa.
\end{split}
\end{equation}
To bound $I_2$ we write $[U,V] = [U,f(X)/\kappa]\cup [f(X)/\kappa,V]$, we apply partial integration to the first range, then we bound $\sin (\kappa t)/t$ by $1/t$ in the second range and apply partial integration to the term involving $\Phi(t)$. Hence we obtain
\begin{equation}
\label{3-5}
\begin{split}
I_2 &= \Odig{\frac{V\log V}{V^2}} + o(\kappa^2 U\log U) + \Odig{\kappa\int_U^{f(X)/\kappa} |r(t)| \min(\kappa t^{-1};t^{-2}) \dx t} 
\\&
\hskip1cm+ 
\Odig{ \int_{f(X)/\kappa}^{V} \frac{\Phi(t)+ \log t+1}{t^{2}} \dx t} 
 \\
&= \Odig{\kappa\int_U^{f(X)/\kappa}   |r(t)| \min(\kappa t^{-1};t^{-2}) \dx t} +  o\big(\kappa\log(1/\kappa)\big) \\
&= o\big(\kappa\log(1/\kappa)\big)
\end{split}
\end{equation}
provided $f(X)\to\infty$ slowly enough. Finally, using the bound for $\Phi(t)$ in Lemma 3 we obtain
\begin{equation}
\label{3-6}
I_3 \ll \frac{1}{\tau} \int_V^\infty \frac{\log^2t}{t^2} \dx t \ll \kappa,
\end{equation}
and Lemma 8 follows from \eqref{3-1} and \eqref{3-4}-\eqref{3-6}. \fine

\medskip
The next lemma is a $\tau$-uniform version of Lemma 10 of \cite{Go-Mo/1987}. For $s\in\CC$ we write
\begin{equation}
\label{3-7}
c(\theta,s) =  ((1+\theta)^s-1)/s.
\end{equation}

\smallskip
{\bf Lemma 9.} {\sl Let $0 \leq \theta\leq \tau \leq 1$ and $Z\geq 1/\theta$, and let $\Phi(X, t, \tau)$ be as in \eqref{2-1}. Then}
\begin{equation}
\label{3-8}
\begin{split}
\int_{-\infty}^{\infty}  \vert c(\theta,it) \vert^{2} \Phi(X, t, \tau)^2 \dx t = \int_{-\infty}^{\infty} &\Bigl \vert \sum_{\vert\gamma\vert \leq Z} \!\! \frac{c(\theta,1/2+i\gamma)\ X^{i\gamma}}{1+\tau^{2}(t-\gamma)^{2}}\Bigr \vert^{2} \dx t \\
&+ \Odig{\frac{\theta^2}{\tau^3}  \log^{4}(2/\theta)+ \frac{\log^{4}(2Z)}{\tau^{2}Z}}.
\end{split}
\end{equation}

\medskip
{\it Proof.} Let $I$ denote the integral on the left hand side of \eqref{3-8}, and $J$, $K$ be the corresponding integrals with $c(\theta,it)$ replaced by $c(\theta,1/2+it)$ and $c(\theta,1/2+i\gamma)$, respectively (clearly, $c(\theta,1/2+i\gamma)$ is inside the sum defining $\Phi(X,t,\tau)$, see \eqref{2-1}). Write $J = \int_\R \vert A \vert ^{2}$ and $K = \int_\R \vert B \vert ^{2}$. In the strips $\vert \Re(s) \vert \leq 1/\theta$ and $\vert \Re(s) \vert \leq 1/(2\theta)$ 
we have, respectively,
\begin{equation}
\label{3-9}
c(\theta,s)\ll \min(\theta;\vert s \vert ^{-1})
\quad\textrm{and}\quad
\frac{\partial c(\theta,s)}{\partial s} \ll \theta \min(\theta;\vert s \vert ^{-1}).
\end{equation}
Since $|\alpha|^2-|\beta|^2\leq |\alpha-\beta|(|\alpha|+|\beta|)$ for $\alpha,\beta\in\CC$, a direct application of \eqref{3-9} and Lemma 3, recalling that $\tau \ge \theta$, gives
\begin{equation}
\label{3-10}
I-J \ll \frac{\theta}{\tau^2} \int_{-\infty}^\infty \min\Bigl(\theta;\frac{1}{|t|}\Bigr)^2 \log^{2} \max (\vert t\vert +2; 2/\tau) \dx t \ll \frac{\theta^{2}}{\tau^{2}} \log^{2}(2/\theta).
\end{equation} 

\smallskip
By a similar but more complicated argument, coupled with Lemma 2, we also obtain
\begin{equation}
\label{3-11}
 J-K \ll \frac{\theta^2}{\tau^{3}}  \log^{4}(2/\theta).
\end{equation} 
Indeed, \eqref{3-9} and Lemma 3 immediately give
\[
 A \ll \frac{1}{\tau} \min\Bigl( \theta; \frac{1}{\vert t \vert} \Bigr) \log \max ( \vert t\vert +2;2/\tau)
\] 
and
\[
B \ll \min\Bigl\{ \frac{\theta}{\tau} \log\max (\vert t\vert +2;2/\tau) ; \sum_{\gamma}
 \frac{1}{\vert \gamma\vert} \frac{1}{1+\tau^{2}(t-\gamma)^{2}} \Bigr\}.
 \]
The part with $|\gamma|>2|t|$ of the above sum is $\ll \frac{1}{\tau \vert t \vert} \log \max (\vert t\vert +2;2/\tau)$, while
\[
\sum_{\vert\gamma\vert \leq 2\vert t \vert} \frac{1}{\vert \gamma\vert} \frac{1}{1+\tau^{2}(t-\gamma)^{2}}  \ll  \log (\vert t \vert +2) \int_{1}^{\vert t \vert - 1} \frac{\dx u}{(\vert t \vert - u) (1+\tau^{2} u^2)}  +  \frac{1}{\tau \vert t \vert} \log \max (\vert t\vert +2;2/\tau).
\]
Writing
\[
 \frac{1}{(\vert t \vert - u) (1+\tau^{2} u^2)} = \frac{a}{|t|-u} + \frac{bu+c}{1+\tau^{2} u^2},
\] 
 computing $a,b,c$ and then computing the integral we get
 \begin{equation}
 \label{3-12}
\sum_{\vert\gamma\vert \leq 2\vert t \vert} \frac{1}{\vert \gamma\vert} \frac{1}{1+\tau^{2}(t-\gamma)^{2}}  \ll \log (\vert t \vert +2) \frac{\tau |t| + \log (|t|+2)}{1+\tau^2t^2} \ll \frac{1}{\tau \vert t \vert} \log^2 \max (\vert t\vert +2;2/\tau),
\end{equation}
provided $|t|\geq 1/\theta \ (\geq 1/\tau)$. Collecting the above bounds we obtain
 \begin{equation}
 \label{3-13}
 A+B \ll \frac{1}{\tau} \min\Bigl( \theta; \frac{1}{\vert t \vert} \Bigr) \log^{2} \max  (\vert t\vert +2;2/\tau).
\end{equation}

\smallskip
Turning to the estimate of $A-B$, if $|t|>1/\theta$ we use the mean-value theorem and both bounds in \eqref{3-9} when $|\gamma|\leq2|t|$ and $t\gamma>0$, otherwise only the first one. Therefore
\[
\begin{split}
A-B &\ll \theta  \sum_{|\gamma|\leq2|t|, t\gamma>0} \frac{|t-\gamma|}{|\gamma|(1+\tau^2(t-\gamma)^2)} +\theta \sum_{|\gamma|\leq 2|t|, t\gamma<0}  \frac{1}{1+\tau^2(t-\gamma)^2} \\
& \hskip.5cm+  \theta \sum_{|\gamma|> 2|t|} \frac{1}{1+\tau^2(t-\gamma)^2} = \Sigma_1 + \Sigma_2 + \Sigma_3,
\end{split}
\]
say. Arguing similarly as in the previous case we obtain
\[
\Sigma_1 \ll \frac{\theta}{\tau^2}\log(|t|+2) \int_1^{|t|-1} \Bigl\{\frac{1}{(|t|+u)u} + \frac{1}{(|t|-u)u} \Bigr\} \dx u \ll \frac{\theta}{\tau^2|t|}\log^2(|t|+2),
\]
and using Lemma 2 we also get
\[
\Sigma_2 + \Sigma_3  \ll \theta \sum_{n\geq |t|} \frac{\log n}{1+\tau^2n^2} \ll \frac{\theta}{\tau^2|t|}\log(|t|+2).
\]
If $|t|\leq 1/\theta$, in a similar way we have, recalling that $\tau\geq \theta$,
\[
A-B \ll \theta^2 \sum_{|\gamma|\leq 2/\theta} \frac{|t-\gamma|}{1+\tau^2(t-\gamma)^2} + \theta\sum_{|\gamma|>2/\theta} \frac{1}{1+\tau^2(t-\gamma)^2} \ll \frac{\theta^2}{\tau^2} \log^2(2/\theta).
\]
Collecting the above estimates we get
 \begin{equation}
 \label{3-14}
 A-B \ll \frac{\theta}{\tau^2} \min\Bigl( \theta; \frac{1}{\vert t \vert} \Bigr) \log^{2} \max (\vert t\vert +2;2/\theta),
\end{equation}
and \eqref{3-11} follows immediately from \eqref{3-13} and \eqref{3-14}.

\smallskip
Finally let $L=\int_\R \vert C \vert ^{2}$ be the integral on the right hand side of \eqref{3-8}. Then arguing as before we obtain
 \begin{equation}
 \label{3-15}
 B+C \ll \frac{1}{\tau} \min\Bigl(\theta; \frac{1}{\vert t \vert}\Bigr) \log^2 \max (\vert t\vert +2;2/\tau).
\end{equation}
Moreover, we have (recall that $Z\geq 1/\theta$)
 \begin{equation}
 \label{3-16}
B-C \ll \frac{\log\max (\vert t \vert+2;2/\tau)}{\tau}
\begin{cases}
1/Z &  \textrm{if}\
 Z> 2 \vert t \vert \\
\frac{\log\max (\vert t \vert+2;2/\tau)}{\vert t \vert}&
 \textrm{otherwise}.
\end{cases} 
\end{equation}
Indeed by \eqref{3-9} we have
\[
B-C \ll \sum_{|\gamma|>Z} \frac{1}{|\gamma|(1+\tau^2(t-\gamma)^2)},
\]
hence if $Z>2|t|$ by Lemma 3 we get
\[
B-C \ll \frac{1}{\tau Z} \log\max (\vert t \vert+2; 2/\tau),
\]
while if $Z\leq 2|t|$ we have
\[
B-C \ll \big(\sum_{Z<|\gamma|\leq 2|t|} + \sum_{|\gamma|>2|t|}\big) \frac{1}{|\gamma|(1+\tau^2(t-\gamma)^2)} \ll \frac{1}{\tau|t|} \log^2\max(|t|+2;2/\tau)
\]
thanks to \eqref{3-12}. Therefore, \eqref{3-15} and \eqref{3-16} give
\begin{equation}
\label{3-17}
  K-L \ll \frac{1}{\tau^{2}Z} \log^{4}(2Z),
\end{equation} 
and the lemma follows from \eqref{3-10}, \eqref{3-11} and \eqref{3-17}. \fine

\medskip
The last lemma gives a mean-square estimate for the error term in the explicit formula. Writing
\[
U(X,\tau,\theta) = \int_X^{X(1+\tau)}\Bigl\vert\sum _{\vert \gamma \vert \leq Z}c(\theta,\rho) x^\rho\Bigr\vert^2 \dx x,
\]
and arguing as in Lemma 6 of \cite{L-P-Z/2012} we have

\medskip
{\bf Lemma 10.} {\sl For $Z\geq X\log^{2}X$, $\theta\in(0,1]$ and $\tau\in(0,1]$ we have}
\[
J(X, \tau, \theta)=U(X,\tau,\theta)+\Odig{\frac{\tau X^2\log^{2}(XZ)}{Z}}+\Odig{\frac{U(X,\tau,\theta)^{1/2}\tau^{1/2}X\log (XZ)}{Z^{1/2}}}.
\]

\medskip
{\bf Proof of Theorem 2}. Writing $\kappa=\frac{1}{2} \log(1+\theta)$ we have $\kappa = \theta/2 + \Odi{\theta^{2}}$ and $\log (1/\kappa) = \log(1/\theta)+\log 2 + \Odi{\theta}$, uniformly for $1/X \leq \theta \leq X^{-\eps}$ and $\theta^{1/2-\eps}\leq \tau \leq 1$. The same uniformity ranges for $\theta$ and $\tau$ hold throughout the proof of Theorem 2; we shall not repeat it further. Recalling \eqref{3-7} we have (see p.198 of \cite{Go-Mo/1987})
\[
\vert c(\theta,it) \vert^{2} = 4 \Bigl( \frac{\sin \kappa t}{t} \Bigr)^{2},
\]
hence from Lemma 8 (assuming the Conjecture) and Lemma 9 we obtain (the integrands are even functions of $t$)
\begin{equation}
\label{3-18}
\tau \int_{-\infty}^{\infty} \Bigl \vert \sum_{\vert \gamma \vert \leq Z} \frac{c(\theta,\rho)\ X^{i\gamma}}{1+\tau^{2}(t-\gamma)^{2}} \Bigr \vert^{2} \dx t \sim \frac{\pi}{2} \theta  \log(1/\theta)
\end{equation}
provided 
\begin{equation}
\label{3-19}
Z \geq  \frac{\log^{4} (1/ \theta)}{\tau\theta}.
\end{equation}
But the Fourier transform of 
\[
\sum_{\vert \gamma \vert \leq Z}  \frac{c(\theta, \rho) X^{i\gamma}}{1+\tau^{2}(t-\gamma)^{2}} \qquad \text{is} \qquad
\frac{\pi}{\tau} \sum_{\vert \gamma \vert \leq Z} c(\theta,\rho) X^{i\gamma} e(-\gamma u) e^{-2\pi \vert u \vert/\tau}
\]
(see again p.198 of \cite{Go-Mo/1987}) and hence, by Plancherel's theorem, \eqref{3-18} becomes 
\[
\int_{-\infty}^{\infty} \Bigl \vert \sum_{\vert \gamma \vert \leq Z} c(\theta,\rho) X^{i\gamma} e(-\gamma u) \Bigr \vert^{2} e^{-4\pi \vert u \vert/\tau} \dx u \sim \frac{1}{2\pi} \tau \theta \log(1/\theta).
\]
Therefore by the substitutions $X= e^{\tau Y}$ and $-2\pi u = \tau y$ we have
\begin{equation}
\label{3-20}
\int_{-\infty}^{\infty} \Bigl \vert  \sum_{\vert \gamma \vert \leq Z} c(\theta,\rho) e^{i\gamma \tau (Y+y)} \Bigr \vert^{2} e^{-2\vert y \vert} \dx y \sim \theta  \log(1/\theta).
\end{equation}

We normalize the integrand in \eqref{3-20} dividing by $\theta \log (1/ \theta)$ and then apply  Lemma 7 to the normalized integrand, with $1/X \leq \theta \leq X^{-\eps}$ and $Z\geq  (\tau\theta)^{-1}\log^{4} (1/ \theta)$. Multiplying again by $\theta\log (1/ \theta)$ we obtain 
\[
 \int_0^{\frac{\log (1+\tau)}{\tau}} \Bigl \vert \sum_{\vert \gamma \vert \leq Z} c(\theta,\rho) e^{i\gamma\tau(Y+y)} \Bigr \vert^{2} \, e^{2y\tau} \dx y \sim \Bigl(1+\frac{\tau}{2}\Bigr)\, \theta  \log(1/\theta),
\]
and substituting $x=\exp(\tau(Y+y))$ we get
\begin{equation}
\label{3-21}
 \int_{X}^{X(1+\tau)}  \Bigl \vert  \sum_{\vert \gamma \vert \leq Z} c(\theta,\rho) x^{\rho} \Bigr \vert^{2} \dx x \sim \Bigl(1+\frac{\tau}{2}\Bigr) \, \tau \theta X^{2}  \log(1/\theta).
\end{equation}
This corresponds to the third displayed equation on p.199 of \cite{Go-Mo/1987}. Theorem 2 follows combining \eqref{3-21} and Lemma 10 with the choice $Z =  X \tau^{-1}\log^{4} X$, which is admissible for $1/X \leq \theta \leq X^{-\eps}$ in view of \eqref{3-19}.

%%%%%%%%%%%%%%%%%%%%%%%%%%%%%%%%%%%%%%%%%%%%%
%%%%%%%%%%%%%%%%%%%%%%%%%%%%%%%%%%%%%%%%%%%%

\vskip1cm
%\newpage

\vskip1cm 
\noindent
Alessandro Languasco, Dipartimento di Matematica, Universit\`a
di Padova, Via Trieste 63, 35121 Padova, Italy. \url{languasco@math.unipd.it}

\smallskip
\noindent
Alberto Perelli, Dipartimento di Matematica, Universit\`a di Genova, via
Dodecaneso 35, 16146 Genova, Italy. \url{perelli@dima.unige.it}

\smallskip
\noindent
Alessandro Zaccagnini, Dipartimento di Matematica e Informatica, Universit\`a di Parma, Parco
Area delle Scienze 53/a, 43124 Parma, Italy. \url{alessandro.zaccagnini@unipr.it}
\end{document}